\pdfoutput=1
\RequirePackage{ifpdf}
\ifpdf 
\documentclass[pdftex]{sigma}
\else
\documentclass{sigma}
\fi

\begin{document}

\newcommand{\arXivNumber}{1308.1141}

\allowdisplaybreaks

\renewcommand{\thefootnote}{$\star$}

\renewcommand{\PaperNumber}{094}

\FirstPageHeading

\ShortArticleName{$\mathcal{A}=\mathcal{U}$ for Locally Acyclic Cluster Algebras}

\ArticleName{$\boldsymbol{\mathcal{A}=\mathcal{U}}$ for Locally Acyclic Cluster Algebras\footnote{This paper is
a~contribution to the Special Issue on New Directions in Lie Theory.
The full collection is available at
\href{http://www.emis.de/journals/SIGMA/LieTheory2014.html}{http://www.emis.de/journals/SIGMA/LieTheory2014.html}}}

\Author{Greg MULLER}

\AuthorNameForHeading{G.~Muller}

\Address{Department of Mathematics, Louisiana State University, USA}
\Email{\href{mailto:morilac@umich.edu}{morilac@umich.edu}}
\URLaddress{\url{http://www-personal.umich.edu/~morilac/}}

\ArticleDates{Received May 16, 2014, in f\/inal form August 25, 2014; Published online September 03, 2014}

\Abstract{This note presents a~self-contained proof that acyclic and locally acyclic cluster algebras coincide with
their upper cluster algebras.}


\Classification{13F60; 13B30}

\renewcommand{\thefootnote}{\arabic{footnote}}
\setcounter{footnote}{0}

\section{Introduction}

\emph{Cluster algebras} are commutative domains with distinguished generators (cluster variables) and certain
combinatorial identities between them (mutation relations).
They were introduced in~\cite{FZ02} to study dual canonical bases for Lie groups, and have since been discovered in
a~remarkable range of applications (see~\cite{FZ03b,GSV10,Kel11b} for surveys).

\subsection{Upper cluster algebras}

Each cluster algebra $\mathcal{A}$ also determines an \emph{upper cluster algebra} $\mathcal{U}$, with
$\mathcal{A}\subseteq \mathcal{U}$.
Def\/ined as an intersection of Laurent rings, the upper cluster algebra is more natural than the cluster algebra from
a~geometric perspective.

Upper cluster algebras were introduced in the seminal paper~\cite{BFZ05}, where the authors sought to show that certain
rings are cluster algebras\footnote{Specif\/ically, coordinate rings of double Bruhat cells.}. They instead proved that
these rings are upper cluster algebras.
Naturally, they also asked when $\mathcal{A}\subseteq\mathcal{U}$ is equality.

To this end, they introduced \emph{acyclic cluster algebras}: a~class of elementary cluster algebras which have since
proven to be particularly easy to work with.
They also def\/ined a~\emph{totally coprime} condition; a~more technical condition that can depend on the coef\/f\/icients of
the cluster algebra.
Using both properties, they were able to trap $\mathcal{A}$ and $\mathcal{U}$ between two simpler algebras (the
\emph{lower bound} and the \emph{upper bound}) and `close the gap' between the bounds to show that
$\mathcal{A}=\mathcal{U}$.

\begin{theorem}[\cite{BFZ05}]
\label{thm: BFZ}
If $\mathcal{A}$ is a~totally coprime, acyclic cluster algebra, then $\mathcal{A}=\mathcal{U}$.
\end{theorem}

\subsection{Locally acyclic cluster algebras}

The current paper's author became interested in studying when $\mathcal{A}=\mathcal{U}$ for cluster algebras coming from
marked surfaces\footnote{$\mathcal{A}=\mathcal{U}$ for marked surfaces was needed to prove that \emph{skein algebras}
are cluster algebras (see~\cite{MulSk}).}. Fortuitously, the geometric techniques developed for that problem (cluster
localization and covers) generalized beyond marked surfaces.
In~\cite{MulLA}, \emph{locally acyclic cluster algebras} were def\/ined to capture those cluster algebras for which these
techniques could show $\mathcal{A}=\mathcal{U}$.

Unfortunately, the proof that $\mathcal{A}=\mathcal{U}$ for locally acyclic cluster algebras which appears
in~\cite[Theorem~4.1]{MulLA} depends an incorrectly stated version of Theorem~\ref{thm: BFZ}; specif\/ically, the totally
coprime hypothesis was omitted.
This has led to some confusion about whether $\mathcal{A}=\mathcal{U}$ for locally acyclic cluster
algebras\footnote{The author apologizes for any confusion this may have caused.}. Thankfully, the techniques
of~\cite{MulLA} can be used to show that $\mathcal{A}=\mathcal{U}$ for locally acyclic cluster algebras without assuming
totally coprime or using Theorem~\ref{thm: BFZ}\footnote{Though, some lemmas here are inf\/luenced by lemmas
in~\cite{BFZ05}; cf.~Proposition~\ref{prop: A=Uiso}.}; this was alluded to in~\cite[Remark~6.7]{MulLA} but not shown.
Since acyclic cluster algebras are locally acyclic, this shows that the totally coprime assumption may be removed from
Theorem~\ref{thm: BFZ}.

\begin{theorem}
If $\mathcal{A}$ is an acyclic or locally acyclic cluster algebra, then $\mathcal{A}=\mathcal{U}$.
\end{theorem}

The purpose of this note is to present an elementary proof of this fact, which assumes nothing except the Laurent
phenomenon.
This not only resolves the dependency error in~\cite{MulLA}, but serves as a~short and self-contained introduction to
the techniques and ef\/fectiveness of cluster localization and locally acyclic cluster algebras.
The reader might f\/ind this a~more straight-forward and accessible motivation for the study of locally acyclic cluster
algebras than~\cite{MulLA}, which relies heavily on geometric intuition and techniques.

We also work in the setting of cluster algebras with \emph{normalized coefficients}, which are a~bit more general than
the cluster algebras with \emph{geometric coefficients} studied in~\cite{MulLA}.

\section{Cluster algebra recollections}

We recall the def\/inition of cluster algebras with normalized coef\/f\/icients.
This generalizes cluster algebras with geometric coef\/f\/icients; see~\cite{FZ07} for the appropriate correspondence.

\subsection{Seeds and mutation}

Let $\mathbb{P}$ be a~\emph{semifield}: a~torsion-free abelian group (written multiplicatively) equipped with an
\emph{auxiliary addition} $\oplus$ which is commutative, associative and distributive over multiplication.
Its integral group ring $\mathbb{ZP}$ will be the \emph{coefficient ring} of the cluster algebra.

Let $\mathcal{F}$ be a~f\/ield which contains $\mathbb{ZP}$.
A~\textit{seed} of rank~$n$ in $\mathcal{F}$ is a~triple $(\mathbf{x},\mathbf{y},\mathsf{B})$ consisting of three parts:
\begin{itemize}\itemsep=0pt
\item the \textit{cluster}: $\mathbf{x}=\{x_1,x_2,\dots ,x_n\}$ is an~$n$-tuple in $\mathcal{F}$ which freely generates
$\mathcal{F}$ as a~f\/ield over the fraction f\/ield of $\mathbb{ZP}$,
\item the \textit{coefficients}: $\mathbf{y}=\{y_1,y_2,\dots ,y_n\}$ is an~$n$-tuple in $\mathbb{P}$, and
\item the \textit{exchange matrix}: $\mathsf{B}$ is an integral, skew-symmetrizable\footnote{\emph{Skew-symmtrizable}
means there is a~diagonal matrix~$\mathsf{D}$ such that $\mathsf{D}B$ is skew-symmetric.} $n\times n$ matrix.
\end{itemize}
A~seed $(\mathbf{x},\mathbf{y},\mathsf{B})$ may be \textit{mutated} at an index $1\leq k\leq n$, to produce a~new seed
$(\mu_k(\mathbf{x}),\mu_k(\mathbf{y})$, $\mu_k(\mathsf{B}))$ def\/ined as follows:
\begin{itemize}\itemsep=0pt
\item $\mu_k(\mathbf{x}):=\{x_1,x_2,\dots ,x_{k-1},x_k',x_{k+1},\dots ,x_n\}$, where
\begin{gather*}
x_k':=\frac{y_k\prod x_j^{[\mathsf{B}_{jk}]_+} + \prod x_j^{[-\mathsf{B}_{jk}]_+}}{(y_k\oplus 1)x_k}
\end{gather*}
(here, $[\mathsf{B}_{jk}]_+$ denotes $\max(\mathsf{B}_{jk},0)$),
\item $\mu_k(\mathbf{y}):=\{y_1',y_2',\dots ,y_n'\}$, where
\begin{gather*}
y_j':=
\begin{cases}
y_j^{-1} & \text{if}\quad k=j,
\\
y_jy_k^{[\mathsf{B}_{kj}]_+}(y_k\oplus 1)^{-\mathsf{B}_{kj}} & \text{if}\quad  k\neq j,
\end{cases}
\end{gather*}
\item $\mu_k(\mathsf{B})$ is def\/ined entry-wise~by
\begin{gather*}
\mu_k(\mathsf{B})_{ij} =
\begin{cases}
-\mathsf{B}_{ij} & \text{if}\quad k=i \quad\text{or}\quad k=j,
\\
\mathsf{B}_{ij}+ \frac{1}{2}(|\mathsf{B}_{ik}|\mathsf{B}_{kj}+\mathsf{B}_{ik}|\mathsf{B}_{kj}|) & \text{otherwise},
\end{cases}
\end{gather*}
\end{itemize}

Mutating at the same index twice in a~row returns to the original seed.
Permuting the indices $\{1,2,\dots ,n\}$ induces a~new seed in the obvious way.
Two seeds are \textit{mutation-equivalent} if they are related by a~sequence of mutations and permutations.

\subsection{Cluster algebras}

Given a~seed $(\mathbf{x},\mathbf{y},\mathsf{B})$, the union of all the clusters which appear in mutation-equivalent
seeds def\/ines a~set of \textit{cluster variables} in the embedding f\/ield $\mathcal{F}$.
The \textit{cluster algebra} $\mathcal{A}(\mathbf{x},\mathbf{y},\mathsf{B})$ determined~by
$(\mathbf{x},\mathbf{y},\mathsf{B})$ is the unital subring of $\mathcal{F}$ generated by $\mathbb{ZP}$ and the cluster
variables.
The cluster algebra only depends on the mutation-equivalence class of the initial seed, and so the initial seed
$(\mathbf{x},\mathbf{y},\mathsf{B})$ will often be omitted from the notation.

We say $\mathbf{x}=\{x_1,x_2,\dots ,x_n\}$ is a~\emph{cluster} in $\mathcal{A}$ if $\mathbf{x}$ is in some seed
$(\mathbf{x},\mathbf{y},\mathsf{B})$ of $\mathcal{A}$.
A~fundamental property of cluster algebras is the Laurent phenonemon, which states that the localization of
$\mathcal{A}$ at a~cluster $\mathbf{x}$ is the ring of Laurent polynomials in $\mathbf{x}$ over $\mathbb{ZP}$
\begin{gather*}
\mathcal{A}\hookrightarrow \mathcal{A}\big[x_1^{-1},x_2^{-1},\dots ,x_n^{-1}\big] =
\mathbb{ZP}\big[x_1^{\pm1},x_2^{\pm1},\dots ,x_n^{\pm1}\big].
\end{gather*}
Every cluster in $\mathcal{A}$ def\/ines such an inclusion.
Def\/ine the \textit{upper cluster algebra} $\mathcal{U}$ of $\mathcal{A}$ to be the intersection of each of these Laurent
rings, taken inside the common fraction f\/ield $\mathcal{F}$
\begin{gather*}
\mathcal{U}:= \bigcap_{\text{clusters}\; \mathbf{x}\subset \mathcal{A}}
\mathbb{ZP}\big[x_1^{\pm1},x_2^{\pm1},\dots ,x_n^{\pm1}\big]\subset \mathcal{F}.
\end{gather*}
By the Laurent phenomenon, there is an embedding $\mathcal{A}\subseteq \mathcal{U}$.
This inclusion is not always equality (see~\cite[Proposition 1.26]{BFZ05}), but it is an equality in many examples of
cluster algebras, and it is hoped to be an equality in many more important examples.

\begin{remark}
An explanation of the geometric signif\/icance of the upper cluster algebra can be found in~\cite[Section 3.2]{MM14}.
\end{remark}

\section{Cluster localization and covers}

This section reviews the techniques of cluster localization and covers, def\/ined in~\cite{MulLA}.

\subsection{Freezing}
`Freezing' a~cluster variable in a~seed promotes one or more cluster variables to the coef\/f\/icient ring.
The name is motivated by the case of geometric seeds, in which the generators of $\mathbb{P}$ are regarded as `frozen'
cluster variables.

Let $(\mathbf{x},\mathbf{y},\mathsf{B})$ be a~seed of rank~$n$ over a~coef\/f\/icient ring $\mathbb{ZP}$.
Def\/ine the \textit{freezing} of $(\mathbf{x},\mathbf{y},\mathsf{B})$ at $x_n$ to be the seed
$(\mathbf{x}^\dagger,\mathbf{y}^\dagger,\mathsf{B}^\dagger)$ of rank $n-1$ def\/ined as follows.
\begin{itemize}\itemsep=0pt
\item $\mathbb{P}^\dagger:=\mathbb{P}\times \mathbb{Z}$, the direct product of $\mathbb{P}$ with a~free abelian group
$\mathbb{Z}$, whose generator will be denoted $x_n$.
The auxilliary addition is def\/ined~by
\begin{gather*}
\big(p_1x_n^a\big)\oplus \big(p_2x_n^b\big):= (p_1\oplus p_2)x_n^{\min(a,b)}.
\end{gather*}
It follows that $\mathbb{ZP}^\dagger\simeq \mathbb{ZP}[x_n^{\pm1}]$ as rings.
\item There is an obvious isomorphism from the f\/ield $\mathcal{F}\simeq \mathbb{Q}(\mathbb{P},x_1,x_2,\dots ,x_n)$ to the
f\/ield $\mathcal{F}^\dagger\simeq \mathbb{Q}(\mathbb{P}^\dagger,x_1,x_2,\dots ,x_{n-1})$.
Identify $x_i$ with its image under this isomorphism, and let $\mathbf{x}^\dagger:=\{x_1,x_2,\dots ,x_{n-1}\}$.
\item $y_i^\dagger:= y_ix_n^{B_{ni}}$ and $\mathbf{y}^\dagger:=\big\{y_1^\dagger,y_2^\dagger,\dots ,y_{n-1}^\dagger\big\}$.
\item Let $\mathsf{B}^\dagger$ be the submatrix of $\mathsf{B}$ obtained by deleting the~$n$th row and column.
\end{itemize}
Freezing an arbitrary variable $x_i\subset \mathbf{x}$ is def\/ined by conjugating the above construction by any
permutation that sends $x_i$ to $x_n$.
The following proposition is straightforward.

\begin{proposition}

Let $(\mu(\mathbf{x}),\mu(\mathbf{y}),\mu(\mathsf{B}))$ be the mutation of $(\mathbf{x},\mathbf{y},\mathsf{B})$ at the
variable $x_i$, and let $(\mathbf{x}^\dagger,\mathbf{y}^\dagger,\mathsf{B}^\dagger)$ be the freezing of
$(\mathbf{x},\mathbf{y},\mathsf{B})$ at $x_j\neq x_i$.
Then the freezing of $(\mu(\mathbf{x}),\mu(\mathbf{y}),\mu(\mathsf{B}))$ at $x_j$ is the mutation at $x_i$ of
$(\mathbf{x}^\dagger,\mathbf{y}^\dagger,\mathsf{B}')$.
\end{proposition}
It follow that the exchange graph\footnote{The \emph{exchange graph} of $\mathcal{A}$ is the possibly inf\/inite graph
with a~vertex for each seed (up to permutation) and an edge for each mutation between seeds.} of
$(\mathbf{x}^\dagger,\mathbf{y}^\dagger,\mathsf{B}^\dagger)$ is the subgraph of the exchange graph of
$(\mathbf{x},\mathbf{y},\mathsf{B})$ which contains the initial seed and avoids mutating at $x_i$.

Given a~subset $S\subset \mathbf{x}$, the \textit{freezing} of $(\mathbf{x},\mathbf{y},\mathsf{B})$ at~$S$ is the
iterated freezing of each variable in~$S$ in any order.
As above, freezing multiple variables~$S$ commutes with mutation away from~$S$.

\subsection{Cluster localization}

The ef\/fect of freezing a~seed on the corresponding cluster algebras and upper cluster algebras is given by the following
sequence of containments.

\begin{lemma}
\label{lemma: containments}
Let $\{x_{i_1},x_{i_2},\dots ,x_{i_k}\}\subset \mathbf{x}$ be a~set of variables in a~seed
$(\mathbf{x},\mathbf{y},\mathsf{B})$, with freezing $(\mathbf{x}^\dagger,\mathbf{y}^\dagger,\mathsf{B}^\dagger)$.
Let $\mathcal{A}$ and $\mathcal{U}$ be the cluster algebra and upper cluster algebra of
$(\mathbf{x},\mathbf{y},\mathsf{B})$ and let $\mathcal{A}^\dagger$ and $\mathcal{U}^\dagger$ be the cluster algebra and
upper cluster algebra of $(\mathbf{x}^\dagger,\mathbf{y}^\dagger,\mathsf{B}^\dagger)$.
Then there are inclusions in $\mathcal{F}$
\begin{gather*}
\mathcal{A}^\dagger \subseteq \mathcal{A}\big[(x_{i_1}x_{i_2}\dots x_{i_k})^{-1}\big] \subseteq \mathcal{U}\big[(x_{i_1}x_{i_2}\dots
x_{i_k})^{-1}\big]\subseteq \mathcal{U}^\dagger.
\end{gather*}
\end{lemma}
\begin{proof}
The coef\/f\/icient ring $\mathbb{ZP}^\dagger=\mathbb{ZP}\big[x_{i_1}^{\pm1},\dots ,x_{i_k}^{\pm1}\big]$, so $\mathbb{ZP}^\dagger
\subset \mathcal{A}\big[(x_{i_1}x_{i_2}\dots x_{i_k})^{-1}\big]$.
Since the cluster variables of $\mathcal{A}$ contain the cluster variables of $\mathcal{A}^\dagger$, this implies the
f\/irst inclusion.

Since the clusters of $\mathcal{A}$ contain the clusters of $\mathcal{A}^\dagger$, there is a~containment
$\mathcal{U}\subset\mathcal{U}^\dagger$.
Furthermore, $(x_{i_1}x_{i_2}\dots x_{i_k})^{-1}\in \mathbb{P}^\dagger\subset \mathcal{U}^\dagger$, which implies the
last inclusion.

The middle inclusion follows from the Laurent phenomenon.
\end{proof}

Each of these inclusions is fairly interesting; but for the purposes at hand, we are most interested in when the f\/irst
inclusion is equality.

\begin{definition}
If $(\mathbf{x}^\dagger,\mathbf{y}^\dagger,\mathsf{B}^\dagger)$ is a~freezing of $(\mathbf{x},\mathbf{y},\mathsf{B})$ at
$\{x_{i_1},x_{i_2},\dots ,x_{i_k}\}$ such that
\begin{gather*}
\mathcal{A}(\mathbf{x}^\dagger,\mathbf{y}^\dagger,\mathsf{B}^\dagger) =
\mathcal{A}(\mathbf{x},\mathbf{y},\mathsf{B})\big[(x_{i_1}x_{i_2}\dots x_{i_k})^{-1}\big],
\end{gather*}
we say $\mathcal{A}(\mathbf{x}^\dagger,\mathbf{y}^\dagger,\mathsf{B}^\dagger)$ is a~\textit{cluster localization} of
$\mathcal{A}(\mathbf{x},\mathbf{y},\mathsf{B})$.
\end{definition}

This is the most natural way a~cluster structure can descend to a~cluster structure on a~localization; hence the name.
Whenever
$\mathcal{A}(\mathbf{x}^\dagger,\mathbf{y}^\dagger,\mathsf{B}^\dagger)=\mathcal{U}(\mathbf{x}^\dagger,\mathbf{y}^\dagger,\mathsf{B}^\dagger)$,
each inclusion in Lemma~\ref{lemma: containments} is an equality, and so
$\mathcal{A}(\mathbf{x}^\dagger,\mathbf{y}^\dagger,\mathsf{B}^\dagger)$ is a~cluster localization of
$\mathcal{A}(\mathbf{x},\mathbf{y},\mathsf{B})$.

Cluster localizations are transitive; if $\mathcal{A}^\dagger$ is a~cluster localization of $\mathcal{A}$, and
$\mathcal{A}^{\ddagger}$ is a~cluster localization of $\mathcal{A}^\dagger$, then $\mathcal{A}^{\ddagger}$ is a~cluster
localization of $\mathcal{A}$.

\begin{example}
The extreme case $S=\mathbf{x}$ amounts to freezing every variable in a~seed.
The freezing $(\mathbf{x}^\dagger,\mathbf{y}^\dagger,\mathsf{B}^\dagger)=(\varnothing,\varnothing,\varnothing)$ is rank zero
and so
\begin{gather*}
\mathcal{A}\big(\mathbf{x}^\dagger,\mathbf{y}^\dagger,\mathsf{B}^\dagger\big) = \mathbb{ZP}^\dagger \simeq
\mathbb{ZP}\big[x_1^{\pm1},x_2^{\pm1},\dots ,x_n^{\pm1}\big] =
\mathcal{U}\big(\mathbf{x}^\dagger,\mathbf{y}^\dagger,\mathsf{B}^\dagger\big).
\end{gather*}
This is the inclusion coming from the Laurent phenomenon.
\end{example}

\subsection{Covers}

For a~cluster algebra $\mathcal{A}$, a~set $\{\mathcal{A}_i\}_{i\in I}$ of cluster localizations of $\mathcal{A}$ (not
necessarily localizations from the same seed) is called a~\textit{cover} if, for every prime ideal~$P$ in $\mathcal{A}$,
there is some $\mathcal{A}_i$ such that $\mathcal{A}_iP\subsetneq \mathcal{A}_i$.

Covers are transitive; that is, if $\{\mathcal{A}_i\}_{i\in I}$ is a~cover of $\mathcal{A}$, and
$\{\mathcal{A}_{ij}\}_{j\in J_i}$ is a~cover of $\mathcal{A}_i$, then $\bigcup_{i\in I}\{\mathcal{A}_{ij}\}_{j\in J_i}$
is a~cover of $\mathcal{A}$.
However, because there is no `geometric intersection' for cluster localizations, there is no notion of a~common
ref\/inement of two covers.

The following is a~useful property of covers.

\begin{proposition}
\label{prop: cover}
If $\{\mathcal{A}_i\}_{i\in I}$ is a~cover of $\mathcal{A}$, then $\mathcal{A}=\bigcap_{i\in I} \mathcal{A}_i$.
\end{proposition}
\begin{proof}
Let $a\in \bigcap_{i\in I}\mathcal{A}_i$.
Write each $\mathcal{A}_i$ as $\mathcal{A}[d_i^{-1}]$ for some $d_i\in \mathcal{A}$, so that $a\in \bigcap_{i\in I}
\mathcal{A}[d_i^{-1}]$.
For each~$i$, there is some $n_i$ such that $d_i^{n_i}a\in \mathcal{A}$.
Def\/ine the $\mathcal{A}$-ideal
\begin{gather*}
J:= \{b\in \mathcal{A} \,|\, ba\in \mathcal{A}\}.
\end{gather*}
It follows that $d_i^{n_i}\in J$ and $\mathcal{A}_iJ=\mathcal{A}_i$, for all~$i$.
If~$M$ is a~prime $\mathcal{A}$-ideal which contains~$J$, then $\mathcal{A}_iM=\mathcal{A}_i$ for all~$i$.
Hence, no prime $\mathcal{A}$-ideals contain~$J$, and so $J=\mathcal{A}$.
In particular, $1\in J$, and so $a\in \mathcal{A}$.
\end{proof}

The equality $\mathcal{A}=\mathcal{U}$ may be checked locally in a~cover.
\begin{lemma}
\label{lemma: A=Ucover}
Let $\{\mathcal{A}_i\}_{i\in I}$ be a~cover of $\mathcal{A}$.
If $\mathcal{A}_i=\mathcal{U}_i$ for each $i\in I$, then $\mathcal{A}=\mathcal{U}$.
\end{lemma}

\begin{proof}
Lemma~\ref{lemma: containments} implies that $\mathcal{U}\subset \mathcal{U}_i$ for all~$i$, and so
$\mathcal{U}\subseteq \bigcap_{i\in I} \mathcal{U}_i$.
Then
\begin{gather*}
\mathcal{U}\subseteq \bigcap_{i\in I} \mathcal{U}_i = \bigcap_{i\in I} \mathcal{A}_i =\mathcal{A} \subseteq\mathcal{U}
\end{gather*}
and so $\mathcal{A}=\mathcal{U}$.
\end{proof}

\begin{example}
Consider the example with initial seed and cluster algebra
\begin{gather*}
\left(\{x_1,x_2\},\varnothing,\left[\begin{matrix} 0 & -1\\1 & 0\end{matrix}\right]\right),
\qquad
\mathcal{A} = \mathbb{Z}\left[x_1,x_2,\frac{x_2+1}{x_1},\frac{x_1+x_2+1}{x_1x_2},\frac{x_1+1}{x_2} \right].
\end{gather*}
The two freezings of the initial seed both def\/ine cluster localizations
\begin{gather*}
\mathcal{A}\big[x_1^{-1}\big] = \mathbb{Z}[x_1^{\pm1}] \left[x_2,\frac{x_1+1}{x_2}\right],
\qquad
\mathcal{A}\big[x_2^{-1}\big] = \mathbb{Z}[x_2^{\pm1}] \left[x_1,\frac{x_2+1}{x_1}\right],
\end{gather*}
which collectively def\/ine a~cover of $\mathcal{A}$, because no prime ideal in $\mathcal{A}$ can contain both $x_1$ and
$x_2$.
Proposition~\ref{prop: cover} implies that
\begin{gather*}
\mathbb{Z}\big[x_1^{\pm1}\big] \left[x_2,\frac{x_1+1}{x_2}\right] \cap \mathbb{Z}\big[x_2^{\pm1}\big] \left[x_1,\frac{x_2+1}{x_1}\right]
= \mathbb{Z}\left[x_1,x_2,\frac{x_2+1}{x_1},\frac{x_1+x_2+1}{x_1x_2},\frac{x_1+1}{x_2} \right].
\end{gather*}
Since $\mathcal{U}_i=\mathcal{A}_i$ for both cluster localizations, Lemma~\ref{lemma: A=Ucover} implies that
$\mathcal{U}=\mathcal{A}$.
\end{example}

\section{Acyclic and locally acyclic cluster algebras}

This section shows that acyclic cluster algebras admit covers by \emph{isolated cluster algebras}.

\subsection{Isolated cluster algebras}

A cluster algebra is \textit{isolated} if the exchange matrix of any seed (equivalently, every seed) is the zero
matrix\footnote{For a~seed of geometric type def\/ined by a~quiver, being isolated is equivalent to having no arrows
between mutable vertices; this is the origin of the term `isolated'.}.

\begin{proposition}
\label{prop: A=Uiso}
Let $\mathcal{A}$ be an isolated cluster algebra.
Then $\mathcal{A}=\mathcal{U}$.\footnote{This is a~straightforward generalization of the statement and proof
of~\cite[Lemma 6.2]{BFZ05}.}
\end{proposition}
\begin{proof}
Let $(\mathbf{x},\mathbf{y},\mathsf{B})$ be a~seed for $\mathcal{A}$.
Let $\mathbf{x}=\{x_1,x_2,\dots ,x_n\}$ and let $x_i':=P_ix_i^{-1}$ denote the mutation of $x_i$ in the initial seed.
Since $\mathsf{B}=0$, each of the $P_i\in \mathbb{ZP}$, and mutation does not change $\mathbf{y}$ or $\mathsf{B}$.
Hence, mutating at distinct indices is order-independent.
It follows that $\{x_1,x_2,\dots ,x_n,x_1',x_2',\dots ,x_n'\}$ is the complete set of cluster variables in $\mathcal{A}$, and
the clusters are of the form
\begin{gather*}
\{x_i |i\not\in I\} \bigcup \{x_i' \,|\, i\in I\}
\end{gather*}
for any subset $I\subset \{1,2,\dots ,n\}$.

Choose some $a\in \mathcal{U}$.
Since $a\in \mathbb{ZP}\big[x_1^{\pm1},x_2^{\pm1},\dots ,x_n^{\pm1}\big]$, it can be written as
\begin{gather*}
a = \sum\limits_{\alpha\in \mathbb{Z}^n} \lambda_\alpha x_1^{\alpha_1}x_2^{\alpha_2}\cdots x_n^{\alpha_n},
\qquad
\lambda_\alpha \in \mathbb{ZP}.
\end{gather*}
For any $I\subset \{1,2,\dots ,n\}$, the element~$a$ can be written as
\begin{gather*}
a = \sum\limits_{\alpha\in \mathbb{Z}^n} \gamma_{\alpha,I} \bigg(\prod\limits_{i\not\in I}
x_i^{\alpha_i}\bigg)\bigg(\prod\limits_{i\in I} x_i'^{-\alpha_i}\bigg)   = \sum\limits_{\alpha\in \mathbb{Z}^n}
\gamma_{\alpha,I} \bigg(\prod\limits_{i\in I} P_i^{-\alpha_i} \bigg) x_1^{\alpha_1}x_2^{\alpha_2}\cdots x_n^{\alpha_n}.
\end{gather*}
Since the monomials in $x_1,x_2,\dots ,x_n$ are a~basis for $\mathbb{ZP}\big[x_1^{\pm1},x_2^{\pm1},\dots ,x_n^{\pm1}\big]$, it follows
that $\lambda_\alpha = \gamma_{\alpha,I} \big(\prod\limits_{i\in I} P_i^{-\alpha_i} \big)$ for all~$\alpha$.
Hence,~$a$ can be written as
\begin{gather*}
a = \sum\limits_{\alpha\in \mathbb{Z}} \gamma_\alpha \bigg(\prod\limits_{i \,|\, \alpha_i\geq0}
x_i^{\alpha_i}\bigg)\bigg(\prod\limits_{i \,|\, \alpha_i<0} x_i'^{-\alpha_i}\bigg),
\end{gather*}
where $\gamma_\alpha$ is $\gamma_{\alpha,I}$ for~$I$ the subset where $\alpha_i<0$.
This expression for~$a$ is clearly in $\mathcal{A}$.
\end{proof}

\begin{corollary}
\label{coro: A=ULA}
If $\mathcal{A}$ admits a~cover by isolated cluster algebras, then $\mathcal{A}=\mathcal{U}$.
\end{corollary}
\begin{proof}
This is an immediate corollary of Proposition~\ref{prop: A=Uiso} and Lemma~\ref{lemma: A=Ucover}.
\end{proof}

\subsection{Acyclic cluster algebras}

Isolated cluster algebras are too simple to be interesting by themselves, but the key idea of~\cite{MulLA} is that many
interesting cluster algebras can be covered by isolated cluster algebras.

A cluster algebra is \textit{acyclic} if it has an \emph{acyclic seed}: a~seed with an exchange matrix $\mathsf{B}$,
such that there is no sequence of indices $i_1,i_2,\dots ,i_\ell\in \{1,2,\dots ,n\}$ with $\mathsf{B}_{i_{j+1}i_j}>0$ for
all~$j$ and $i_\ell=i_1$.\footnote{For a~seed of geometric type def\/ined by a~quiver, being acyclic is equivalent to
having no directed cycles; this is the origin of the term `acyclic'.
Not every seed of an acyclic cluster algebra is an acyclic seed.} Acyclic cluster algebras form a~class of well-behaved
cluster algebras which contain many notable examples, including all f\/inite-type cluster algebras.

\begin{proposition}
If $\mathcal{A}$ is acyclic, then it admits a~cover by isolated cluster algebras.
\end{proposition}
\begin{proof}
The proof will be by induction on the rank~$n$ of $\mathcal{A}$.
If acyclic $\mathcal{A}$ has rank $n\leq1$, then~$\mathcal{A}$ is isolated, and trivially has a~cover by isolated
cluster algebras.

Assume that every acyclic cluster algebra of rank $<n$ admits a~cover by isolated cluster algebras.
Let $\mathcal{A}$ be an acyclic cluster algebra of rank~$n$, and let $(\mathbf{x},\mathbf{y},\mathsf{B})$ be an acyclic
seed.
There must be some index $i\in \{1,2,\dots ,n\}$ which is a~\emph{sink}; that is, $\mathsf{B}_{ji}\leq0$ for all~$j$.
Otherwise, it would be possible to create arbitrarily long sequences of indices $i_1,i_2,\dots ,i_\ell$ with
$\mathsf{B}_{i_{j+1}i_j}>0$; by f\/initeness, at least one such sequence will have $i_\ell=i_1$.

If $\mathcal{A}$ is isolated, then it trivially has a~cover by isolated cluster algebras.
If $\mathcal{A}$ is not isolated, then there must be indices~$i$ and~$j$ such that~$i$ is a~sink and
$\mathsf{B}_{ji}<0$.
The mutation relation at~$i$ is then
\begin{gather*}
x_ix_i' = \frac{y_i}{y_i\oplus 1} + \frac{1}{y_i\oplus 1}\prod\limits_{k\,|\, \mathsf{B}_{ki}<0} x_k^{-\mathsf{B}_{ki}}.
\end{gather*}
Since $\frac{y_i}{y_i\oplus 1}$ is invertible in $\mathbb{ZP}$, this can be rewritten as
\begin{gather*}
1 = \frac{y_i\oplus 1}{y_i} x_i'x_i - y_i^{-1}\prod\limits_{k\,|\, \mathsf{B}_{ki}<0} x_k^{-\mathsf{B}_{ki}}.
\end{gather*}
Since $x_j$ appears in the right-hand product, $x_i$ and $x_j$ generate the trivial $\mathcal{A}$-ideal.

Let $\mathcal{A}_i$ and $\mathcal{A}_j$ denote the freezings of $\mathcal{A}$ at the indices~$i$ and~$j$, respectively.
The freezing of an acyclic seed is an acyclic seed, and so $\mathcal{A}_i$ and $\mathcal{A}_j$ are acyclic cluster
algebras of rank $n-1$.
By the inductive hypothesis, they admit covers by isolated cluster algebras, and by Corollary~\ref{coro: A=ULA},
$\mathcal{A}_i=\mathcal{U}_i$ and $\mathcal{A}_j=\mathcal{U}_j$.
By Lemma~\ref{lemma: containments}, these are cluster localizations with $\mathcal{A}_i=\mathcal{A}\big[x_i^{-1}\big]$ and
$\mathcal{A}_j=\mathcal{A}_j\big[x_j^{-1}\big]$.

Let~$P$ be a~prime $\mathcal{A}$-ideal.
Since $x_i$ and $x_j$ generate the trivial $\mathcal{A}$-ideal,~$P$ cannot contain both elements.
If $x_i\not\in P$, then $\mathcal{A}_iP \subsetneq \mathcal{A}_i$.
If $x_j\not\in P$, then $\mathcal{A}_jP\subsetneq \mathcal{A}_j$.
Hence, $\mathcal{A}_i$ and $\mathcal{A}_j$ cover $\mathcal{A}$.
The union of the covers of $\mathcal{A}_i$ and $\mathcal{A}_j$ by isolated cluster algebras def\/ines a~cover
$\mathcal{A}$ by isolated cluster algebras.
\end{proof}

\begin{corollary}
If $\mathcal{A}$ is acyclic, then $\mathcal{A}=\mathcal{U}$.
\end{corollary}

\subsection{Locally acyclic cluster algebras}

Many non-acyclic cluster algebras also admit covers by isolated cluster algebras.
As an example, every triangulable marked surface determines a~cluster algebra (see~\cite{FG06,FST08, GSV05}), and these
cluster algebras are locally acyclic as long as there are at least two marked points on the boundary~\cite[Theorem~10.6]{MulLA}.
Since this class is interesting in its own right, we give it a~name.

\begin{definition}[\cite{MulLA}] A~cluster algebra is \textit{locally acyclic} if it admits a~cover by isolated cluster algebras.
\end{definition}

Since isolated cluster algebras are acyclic, and every acyclic cluster algebra admits a~cover by isolated cluster
algebras, this condition is equivalent to admitting a~cover by acyclic cluster algebras; hence the name.
By Corollary~\ref{coro: A=ULA}, every locally acyclic cluster algebra has $\mathcal{A}=\mathcal{U}$.
More generally, locally acyclic cluster algebras possess any property of isolated cluster algebras which is
`geometrically local'.
For example, they are integrally closed, f\/initely generated and locally a~complete intersection; see~\cite{MulLA} for
details.

\subsection*{Acknowledgements}

The author gratefully acknowledges the referees for suggestions incorporated into the f\/inal version of this note.

\pdfbookmark[1]{References}{ref}
\LastPageEnding

\end{document}